\documentclass[letterpaper,10pt,conference]{ieeeconf}

\IEEEoverridecommandlockouts
\newcommand{\mypara}[1]{\smallskip\noindent{\bf {#1}.}~}

\usepackage{csquotes}
\usepackage{dirtytalk}
\usepackage{algorithm}
\usepackage{amssymb}
\usepackage{graphicx}
\usepackage[cmex10]{amsmath}
\usepackage{array}
\usepackage{mdwtab}
\usepackage{fixltx2e}
\usepackage{graphics}
\usepackage{epsfig}
\usepackage{times}
\usepackage{amsmath}
\usepackage{amssymb}
\usepackage{url}
\usepackage{lscape}
\usepackage{color}
\usepackage{epsfig}
\usepackage{cite}
\usepackage{balance}
\usepackage{dblfloatfix}
\usepackage{hyperref}
\usepackage{caption}
\usepackage{subcaption}

\usepackage{tabularx}
\usepackage[noend]{algpseudocode}

\makeatletter
\newcommand{\multiline}[1]{%
  \begin{tabularx}{\dimexpr\linewidth-\ALG@thistlm}[t]{@{}X@{}}
    #1
  \end{tabularx}
}
\makeatother

\title{\LARGE \bf Implementing Optimization-Based Control Tasks in Cyber\textendash Physical Systems With Limited Computing Capacity}

\author{Mehdi Hosseinzadeh, Bruno Sinopoli, Ilya Kolmanovsky, and Sanjoy Baruah
\thanks{This work has been supported by National Science Foundation (NSF) under grant numbers ECCS-1931738, ECCS-1932530, and CMMI-1904394.}
\thanks{M. Hosseinzadeh and B. Sinopoli are with the Department of Electrical and Systems Engineering, Washington University in St. Louis, St. Louis, MO 63130, USA (email: mehdi.hosseinzadeh@ieee.org; bsinopoli@wustl.edu). I. Kolmanovsky is with the Department of Aerospace Engineering, University of Michigan, Ann Arbor, MI 48109, USA (email: ilya@umich.edu). S. Baruah is with the Department of Computer Science and Engineering, Washington University in St. Louis, St. Louis, MO 63130, USA (email: baruah@wustl.edu).}
}

\begin{document}

\maketitle
\thispagestyle{empty}
\pagestyle{empty}

\section{Introduction and Motivation}\label{sec:introduction}
A common aspect of today's Cyber\textendash Physical Systems (CPSs) is that multiple control tasks may execute in a shared processor. For a CPS performing $n$ optimization-based control tasks, denoted by $\tau_k,~k={1,\cdots,n}$, let $\ell_k$ and $\Delta T_k$ denote the worst-case execution time and sampling period of task $\tau_k$, respectively. Under the Earliest Deadline First (EDF) policy, the tasks are schedulable on a single preemptive processor if the following condition \cite{Buttazzo2011} is satisfied:
\begin{align}\label{eq:schedulabilitycondition}
\sum_{k=1}^n \frac{\ell_k}{\Delta T_k}\leq1.
\end{align}
Optimization-based control tasks make use of online optimization and thus have large execution times; hence their sampling periods must be large as well to satisfy condition \eqref{eq:schedulabilitycondition}. However larger sampling periods may cause worse control performance. 

\mypara{Prior Work} Existing methods to address the above-mentioned issue are: i) control-scheduling co-design (e.g., \cite{Pazzaglia2021}), where parameters of control systems are modified at every sampling instant; ii) pre-computation (e.g., \cite{Alessio2009}), where optimal control inputs are computed offline and stored for run-time use; iii) triggering-based control (e.g., \cite{Wang2021}), where a triggering mechanism invokes control tasks; and iv) fixed-iteration optimization (e.g., \cite{Pherson2020}), where control tasks perform a fixed number of iterations to approximately track the solution. These methods either do not guarantee constraint satisfaction (e.g., \cite{Pazzaglia2021}), or do not consider the variability and unpredictability of the task execution time (e.g., \cite{Alessio2009}) and available computing time (e.g., \cite{Wang2021,Pherson2020}). 

Recently, dynamically embedded controllers have been proposed in the control theory literature (e.g., \cite{Nicotra2019,ROTEC}), where the processor, instead of solving an optimization problem, runs a virtual dynamical system whose trajectory converges to the optimal solution. This type of an approach is also pursued in our work limited and variable computing capacity in implementing optimization-based control tasks in CPSs. 

\mypara{Goal} Drawing inspiration from dynamically embedded controllers, the goal of our work is to develop a robust to early termination optimization approach that can be used to effectively solve onboard optimization problems involved in controlling the system despite the presence of unpredictable, variable, and limited computing capacity.

\section{Proposed Solution\textemdash Robust to Early Termination Optimization Approach}\label{sec:solution}
\mypara{Task Details} Suppose that task $\tau_k,~k\in\{1,\cdots,n\}$ solves the following optimization problem:
\begin{align}\label{eq:OptimizationProblem}
\left\{\begin{array}{cl}
     & \min\limits_x f_0(x) \\
  \text{s.t.}   & f_i(x)\leq0,~i=1,\cdots,m
\end{array}
\right.,
\end{align}
where $f_0:\mathbb{R}^p\rightarrow\mathbb{R}$ is a strongly convex objective function to be minimized over the $p$-variable vector $x$, and $f_i(x)\leq0$ is the $i$-th inequality constraint. Note that the mathematical problems in many existing optimization-based control tasks (e.g., Model Predictive Control (MPC) \cite{Rawlings2017}) are in the form of optimization problem \eqref{eq:OptimizationProblem}.

\mypara{Proposed Method} Consider the following \textit{modified} barrier function associated with the optimization problem \eqref{eq:OptimizationProblem}:
\begin{align}\label{eq:BarrierFunction}
\mathcal{B}(x,\lambda)=f_0(x)-\sum\limits_{i=1}^m\lambda_{i}\log(-\beta (f_{i}(x)+1/\beta)+1),
\end{align}
where $\beta\in\mathbb{R}_{>0}$ is the barrier parameter and $\lambda=[\lambda_{1}~\cdots~\lambda_{m}]^\top\in\mathbb{R}^m_{\geq0}$ is the vector of dual parameters.

We consider the following primal-dual gradient flow:
\begin{subequations}\label{eq:systemvirtual}
\begin{align}
\dot{\hat{x}}(t)=&-\sigma\nabla_{\hat{x}}\mathcal{B}\big(\hat{x}(t),\lambda(t)\big),\\
\dot{\hat{\lambda}}_i(t)=&+\sigma\Big(\nabla_{\hat{\lambda}_i}\mathcal{B}\big((\hat{x}(t),\lambda(t)\big)+\Psi_i(t)\Big),
\end{align}
\end{subequations}
where $\sigma\in\mathbb{R}_{>0}$ is a design parameter and $\Psi_i(t)$ is the projection operator onto the normal cone of $\lambda_i$ (see \cite{Nicotra2019}).

\mypara{Convergence}
Let $\beta$ be sufficiently large. Using the following Lyapunov function:
\begin{align}\label{eq:Lyapunov}
V(\hat{x}(t),\hat{\lambda}(t))=&\frac{1}{2\sigma}\left\Vert \hat{x}(t)-x^\ast\right\Vert^2+\frac{1}{2\sigma}\left\Vert \hat{\lambda}(t)-\lambda^\ast\right\Vert^2,
\end{align}
where $(x^\ast,\lambda^\ast)$ is the pair of the optimal solution of \eqref{eq:OptimizationProblem} and the vector of optimal dual parameters, and according to the fact that the operator $\left[\big(\nabla_{\hat{x}}\mathcal{B}(\hat{x},\hat{\lambda})\big)^\top~-\big(\nabla_{\hat{\lambda}}\mathcal{B}(\hat{x},\hat{\lambda})\big)^\top\right]^\top$ is strongly monotone, it can be shown \cite{ROTEC} that $\big(\hat{x}(t),\hat{\lambda}(t)\big)$ exponentially converges to $\big(x^\ast,\lambda^\ast\big)$ as $t\rightarrow\infty$.

\mypara{Constraint-Handling} Since $\mathcal{B}(\hat{x}(t),\hat{\lambda}(t))\rightarrow\infty$ only if $f_i(\hat{x}(t))\rightarrow0^-$ for one or more $i\in\{1,\cdots,m\}$, Alexandrov's theorem \cite{Alexandrov2005} implies that 
\begin{align}\label{eq:limitB1}
\lim\limits_{f_i(\hat{x}(t))\rightarrow0^-\text{ for one or more }i}\;\frac{d}{dt}\mathcal{B}\big(\hat{x}(t),\hat{\lambda}(t)\big)<0,
\end{align}
which asserts that $\mathcal{B}\big(\hat{x}(t),\hat{\lambda}(t)\big)$ must decrease along the system trajectories when these trajectories are near the boundary. Thus, $f_i\big(\hat{x}(t)\big)<0,~i=1,\cdots,m$ for all $t\geq0$. The significance of this conclusion is that the evolution of \eqref{eq:systemvirtual} can be stopped at any time instant with a guaranteed feasible solution.

\mypara{Implementation}
Although system \eqref{eq:systemvirtual} is continuous-time, the above-mentioned properties (i.e., convergence and constraint-handling) are approximately maintained when system \eqref{eq:systemvirtual} is implemented in discrete time by making use of the difference quotient and with a sufficiently small sampling period. Thus, to solve problem \eqref{eq:OptimizationProblem}, one can run system \eqref{eq:systemvirtual} until the available computation time is exhausted (that may not be known in advance), and the solution is sub-optimal and guaranteed to enforce the constraints whenever the evolution of system \eqref{eq:systemvirtual} is terminated. This allows the designer to implement optimization-based control tasks with a small sampling period (and consequently with a minimum degradation in performance), while maintaining optimality and constraint-handling capabilities. It is noteworthy that warm-starting can improve convergence of system \eqref{eq:systemvirtual}.

\section{Experimental Results}\label{sec:simulation}
The objective of this section is to validate the proposed optimization approach and assess its effectiveness. The experiments are carried out on an Intel(R) Core(TM) i7-7500U CPU 2.70 GHz with 16.00 GB of RAM. We use \texttt{YALMIP} toolbox to implement the optimization computations.

We consider a case study where two control tasks are executed on a single processing unit. Tasks $\tau_1$ and $\tau_2$ implement MPC to control DC motors \#1 and \#2 given in \cite{Roy2021}, respectively, such that the control inputs belong to the interval $[-10,10]$. The desired sampling period for both tasks is 20 ms. However, we observe (from 2000 runs) that $\ell_1=\ell_2\approx150$ ms, implying that the schedulability condition \eqref{eq:schedulabilitycondition} cannot be satisfied with the desired sampling periods. To satisfy \eqref{eq:schedulabilitycondition} we would need, for instance, that the tasks execute every 300 ms.

We consider the following cases: i) MPC with sampling period $\Delta T=20$ [ms], which is desired but unimplementable; ii) MPC with sampling period $\Delta T=300$ [ms] to satisfy condition \eqref{eq:schedulabilitycondition}; and iii) MPC with sampling period $\Delta T=20$ [ms] and the proposed optimization approach implemented with $\sigma=10$ and $\beta=10^5$. For comparison purposes, we use the Integral Square Error (ISE) index which is $\text{ISE}(t)\triangleq\int_0^te(\eta)^2d\eta$, where $e(t)$ is the tracking error and the integration is performed since the start of the experiment at time $0$ till the current time $t$.

\begin{figure}[!t]
\centering
\includegraphics[width=8.5cm]{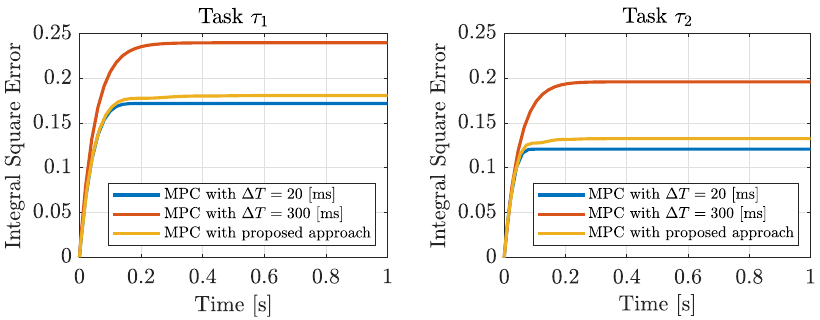}
\caption{Obtained ISEs for considered cases. MPC with the proposed optimization approach has 32\% gain in the obtained ISE for task $\tau_1$ and 48\% gain in the obtained ISE for task $\tau_2$ over the MPC with sampling period $\Delta T=300$ [ms].}
\label{fig:simulation}
\end{figure}

The obtained ISEs for the above-mentioned cases are shown in Fig. \ref{fig:simulation}. As seen in this figure, using a large sampling period degrades the performance by 39\% for task $\tau_1$ and by 62\% for task $\tau_2$. However, MPC with the proposed optimization approach yields better performance by computing a sub-optimal but feasible solution every 20 ms. More precisely, performance degradation with the proposed optimization approach is 4\% for task $\tau_1$ and 9\% for task $\tau_2$.

\section{Conclusion and Future Work}\label{sec:conclusion}
Limited capacity to perform computations is a common constraint in today's CPSs. We proposed an approach to implement optimization-based control tasks despite the variability and unpredictability of available computing time. Experiments were carried out to assess the effectiveness of the proposed approach in satisfying performance requirements and real-time schedulability conditions. Future work will consider how the proposed optimization approach can empower optimization-based control schemes to deal with time-varying constraints.

\bibliographystyle{IEEEtran}
\bibliography{ref}{}

\end{document}